\documentclass{amsart}
\usepackage[russian,english]{babel}
\usepackage{enumitem}
\numberwithin{equation}{section}

\newtheorem{thm}{Theorem}[section]
\newtheorem{lem}[thm]{Lemma}
\newtheorem{cor}[thm]{Corollary}

\newtheorem{defin}[thm]{Definition}
\newtheorem{rem}[thm]{Remark}

\begin{document}
\title{ Time-dependent  identification problem for a fractional Telegraph
equation with the Caputo derivative}

\author{Ravshan Ashurov}
\author{Rajapboy Saparbayev}

\address{Ravshan Ashurov: Institute of Mathematics, Uzbekistan Academy of Science,
Tashkent, Student Town str. 100174.}

\address{ School of Engineering, Central Asian University, 264, Milliy Bog St., 111221, Tashkent, Uzbekistan.}

\email{ashurovr@gmail.com}

\address{Rajapboy Saparbayev: Institute of Mathematics, Uzbekistan Academy of Science,
Tashkent, Student Town str. 100174.}

\email{rajapboy1202@gmail.com}

\small

\title[ Fractional Telegraph
equation ...] {Time-dependent  identification problem for a fractional Telegraph
equation with the Caputo derivative }

\begin{abstract}
This study investigates the inverse problem of determining the right-hand side of a telegraph equation given in a Hilbert space. The main equation under consideration has the form $(D_{t}^{\rho})^{2}u(t)+2\alpha D_{t}^{\rho}u(t)+Au(t)=p( t)q+f(t)$, where $0<t\leq T$, $0<\rho<1$ and $D_{t}^{\rho}$ is the Caputo derivative. The equation contains a self-adjoint positive operator $A$ and a time-varying multiplier $p(t)$ in the source function, which, like the solution of the equation, is unknown. To solve the inverse problem, an additional condition $B[u(t)] = \psi(t)$ is imposed, where $B$ is an arbitrary bounded linear functional. The existence and uniqueness of a solution to the problem are established and stability inequalities are derived. It should be noted that, as far as we know, such an inverse problem for the telegraph equation is considered for the first time. Examples of the operator $A$ and the functional $B$ are discussed.
\end{abstract}
\maketitle

\section{\textbf{Introduction}}

Let $A$ be an arbitrary positive self-adjoint operator defined in some separable Hilbert space $H$ with the domain of definition $D(A)$. Denote the inner product in $H$ by $(\cdot, \cdot)$ and the norm by $||\cdot||$. Let further the inverse operator $A^{-1}$ be compact. Then it is well known, that $A$ has a complete orthonormal system of eigenfunctions
$\{v_k\}$ and a countable set of positive eigenvalues $\lambda_k$ with a unique limit point at 
  $+\infty$. We assume that
the eigenvalues do not decrease as their number increases, i.e. $0<\lambda_1\leq\lambda_2 \cdot\cdot\cdot\rightarrow +\infty$.

For vector functions (or just functions)
$h: \mathbb{R}_+\rightarrow H$ fractional analogs of integrals and derivatives are defined in exactly the same way with scalar functions, and the well-known formulas and
properties are preserved (see, for example, \cite{Liz}). Recall that fractional integrals of the order $\sigma<0$ of the function $h(t)$ defined on $\mathbb{R}_+$ have the form (see, e.g.,~\cite{Pskhu})
\begin{equation}\label{def0}
J_t^\sigma h(t)=\frac{1}{\Gamma
(-\sigma)}\int\limits_0^t\frac{h(\xi)}{(t-\xi)^{\sigma+1}} d\xi,
\quad t>0,
\end{equation}
provided the right-hand side exists. As usual, $\Gamma(\sigma)$ is
Euler's gamma function.The
Caputo fractional derivative of order $\rho\in (0, 1)$ is defined as
$$
D_t^\rho h(t)= J_t^{\rho-1}\frac{d}{dt} h(t).
$$
If $\rho=1$, then the fractional derivative coincides with
the ordinary classical derivative: $D_t h(t)= \frac{d}{dt} h(t)$.

Let $\rho \in (0,1)$ be a fixed number. Consider the following Cauchy problem for the telegraph equation
\begin{equation}\label{prob1}
\begin{cases}
  & (D_{t}^{\rho })^{2}u(t)+2\alpha D_{t}^{\rho }u(t)+Au(t)=p(t)q+f(t), \quad  0<t\le T; \\
 & \underset{t\to 0}{\mathop{\lim }}\,D_{t}^{\rho }u(t)={{\varphi }_{0}},\\
 & u(0)={{\varphi }_{1}}, \\
\end{cases}
\end{equation}
where a part of the source function $p(t)$ is a scalar function, $f(t)$ and $\varphi_{0}, \varphi_{1}, q $ are known elements of $H$.

We note that in most works where the time-fractional telegraph equation is studied, $D_t^{2\rho}$ is taken as the highest derivative instead of $(D_t^{\rho})^2$. In this paper, following paper \cite{RC}, where the motivation for such a choice is also given, the telegraph equation is taken in the form \eqref{prob1}. However, as it was shown in \cite{RC}, these two operators are not the same and the corresponding problems are completely different. As a simple example, we can take the function: 
\begin{equation*}
    u(t)=t^{\rho -1}, \qquad t>0.
\end{equation*}
Simple calculations show that
\[
(D^\rho)^2 u= D^\rho(D^\rho u)\neq D^{2\rho}u.
\]

Problem \eqref{prob1} will be called  \emph{the forward problem} and it was studied in detail in our recent paper \cite{AS}. We will use the results of this work.

The goal of this paper is to find both the solution $u(t)$ and the time-dependent part $p(t)$ of the source function. Since one more unknown is being added, it is natural to add a condition for the unique definition of the new unknown. Following the paper of A. Ashyraliev and H. Al-Khazaima \cite{Ashyr3} (note that in this paper the authors studied the inverse problem for the classical telegraph equation with $\rho=1$) consider the additional condition in a fairly general form:
\begin{equation}\label{ad}
B[u(t)]=\psi(t), \quad 0\leq t \leq T,
\end{equation}
where $B: H\rightarrow \mathbb{R}$ is a given bounded linear functional, and $\psi(t)$ is the given scalar function.

We emphasize, take an additional condition in the form (\ref{ad}), to solve the inverse problem, we come to a rather complicated integro-differential equation involving fractional derivatives. To solve these equations, we use original ideas from the paper of A. Ashyraliev and H. Al-Khazaima \cite{Ashyr3}.

We call the Cauchy problem \eqref{prob1} together with additional condition \eqref{ad}
 \emph{the inverse problem}.

 To the best of our knowledge, the classical telegraph equation first appeared in the work \cite{W.Thomson} of Lord Kelvin in the 19th century, when modeling the movement of an electric current. This is a hyperbolic differential equation with constant coefficients:
\begin{equation}\label{classic}
 u_{tt}+a u_t+b u-cu_{xx}=f(x,t), \quad c>0.   
\end{equation}

Then, in 1876, Oliver Heaviside came up with this equation while simulating the passage of electrical signals in marine telegraph cables (see, e.g.\cite{Lieberstein}, \cite{W.Arendt}). 
In 1958, the same equation was proposed by Cattaneo \cite{C.Cattaneo} to overcome the problem of infinite propagation velocity in heat transfer. Then this equation and its time-fractional versions were arrived at by specialists in modeling various physical processes, such as transmission lines for all frequencies \cite{W.Hayt}, random walks \cite{JBRM}, solar particle transport \cite{F.Effenberger}, oceanic diffusion \cite{A.Okuko}, wave propagation \cite{V.H.Weston}, damped small oscillations, anomalous diffusion and wave processes \cite{L.Boyadjiev,J.Masoliver,Orsingher,Orsingher1}.

In the last two decades, interest in the analysis and applications of fractional, both in time and in space variables, generalizations of telegraph-type equations has increased in the literature. Why use fractional derivatives? The answer lies in the fact that all real physical phenomena that change in time depend not only on the exact moment of time, but also on the previous time, in other words, memory in time. This memory over time is an important aspect of fractional calculus (see \cite{Lorenzo}).
 One of the first works in this direction was the fundamental paper by  R. Cascaval et al. \cite{RC}. In this paper, the highest time derivative is taken in the form $(D_t^{\rho})^2$  (note that in all other papers where the telegraph equation with fractional time derivatives is studied, $D_t^{2\rho}$ is taken) and the rationale for such a choice is given. The main goal of the paper \cite{RC} is to study the asymptotic behavior of the solution $u(t)$ of the homogeneous problem \eqref{prob1} for large $t$. The authors proved the existence of a solution $v(t)$ of the equation $2\alpha D_{t}^{\rho }v(t)+Av(t)=0$ for which the asymptotic relation
\[
u(t)=v(t) + o(v(t)), \quad t\to +\infty,
\]
is valid.

For the telegraph equation with the highest fractional derivative in time $D_t^{2\rho}$ and the elliptic part of the equation in the form $Au(x,t)= u_{xx}(x,t)$, various questions have been considered by a number of researchers (see e.g. \cite{Ashyr3}, \cite{Hashmi},\cite{Doetsch}, \cite{Hosseini}).  Thus, in \cite{Orsingher} the authors found the Fourier transform of the solution of the Cauchy problem. In \cite{Beghin}  (in the case of fractional
derivatives of rational order $\rho=m/n$ with $m < n$), fundamental solutions of the telegraph equation are constructed. The authors of the paper \cite{Orsingher.ZHAO} considered the equation \eqref{classic} with space-fractional derivative $D^\alpha_x u$  of order $1<\alpha<2$ instead of $u_{xx}$  and constructed the Fourier transform of the fundamental solution.
In work \cite{Huang} (see also \cite{YMOstoja}) for the telegraph equation  with $1/2\leq \rho\leq 1$ a fundamental solution of the Cauchy problem in cases $x\in \mathbb{R}$ and $x\in \mathbb{R}_+$ is found using the Fourier-Laplace transforms and their inverse transforms. In the same paper \cite{Huang} for the case of a bounded spatial domain, the formal solution of the boundary value problem is found in the form of a series using the Sine-Laplace transformation method, but convergence of these series is not investigated. Similarly to this paper, the authors of \cite{Ahmad} studied the problem in a bounded spatial domain with the operator $A=(-\triangle) ^{\beta/2}$, $\beta\in ( 0,2] $, and they also found formal analytical solutions under inhomogeneous Dirichlet and Neumann boundary conditions using the method of separation of variables. However, as in \cite{Huang}, the authors did not investigate the convergence and differentiability of these series. 

In all the above works, the Caputo derivative is taken as the fractional derivative. In \cite{RAM.SAX.ORS,MrAjay,VRF} Cauchy problems for telegraph equations are considered, where fractional derivatives with respect to $t$ are taken in the sense of Hilferé or Hadamard, and derivatives with respect to $x$ in the sense of Riesz-Feller. The authors of the paper \cite{RFigueiredo} constracted the Green function for the space-time fractional telegraph equations.

Let us mention one more paper \cite{MAUHafiz}
 (see also the literature there), where the non-linear telegraph equation $$D_t^{2\rho}u-D_x^{2\rho}u+D_t^{\rho}u+ \gamma u +\beta u^3=0$$ 
 is considered. The authors succeeded in constructing a solution of the form $$u=u(y), \, y=k\, \frac{x^\rho}{\rho}-c\,  \frac{t^\rho}{\rho}. $$

A number of specialists have developed numerical algorithms for solving the problem \eqref{prob1} for various  operators $A$. Overview of some
work in this direction is contained in the articles \cite{Momani,Jordan} (see also \cite{NJForda}, where the authors considered the space-time fractional telegraph equations).

In the theory of equations in partial and fractional derivatives, inverse problems are called problems in which, in addition to solving a differential equation, it is required to determine any coefficient (coefficients) of the equation, or its right side, or both coefficients (c) and the right side. Interest in the study of inverse problems for the equations noted above is due to the importance of their applications in various branches of mechanics, seismology, medical tomography, and geophysics. It should be noted that inverse problems are well studied for differential equations of integer order (see, e.g.\cite{Kabanikhin} and the literature therein).

As regards inverse problems for the telegraph equation, the following two papers \cite{LR1,LR2} can be indicated. In work \cite{LR1}, the equation
\[
D_t^\alpha u- r(t) D_t^\beta u+ a^2 (-\triangle)^\gamma u=f(x,t), \quad x\in\Omega,\quad t>0,
\]
with boundary and initial conditions was considered and the existence and uniqueness of the coefficient $r(t)$ under the over-determination condition
\begin{equation}\label{over}
    \int_\Omega u(x,t) \varphi(x) dx = F(t),
\end{equation}
with given functions $\varphi$ and $F$, was proved. A similar inverse problem for an equation whose right-hand side depends nonlinearly on $u$ and $D_t^\beta u$ was studied in \cite{LR2}.

A large number of papers are also devoted to inverse problems of determining the right-hand side of the subdiffusion equations (see, for example,\cite{AOLob,AODif,Kirane,Ruzhansky} and references therein). 

To the best of our knowledge, this paper is the first to study the inverse problem of determining the right-hand side of a time-fractional telegraph equation. It should be specially emphasized that in order to find the additional unknown, we took the over-determination condition of a rather general form. In particular, as the functional $B$ in condition \eqref{ad}, one can take the integral in \eqref{over}.

The work consists of 5 sections and Conclusion. In the next section, we formulate the main results of the work.
Section 3 presents some properties of the Mittag-Leffler functions. To solve the inverse problem, results on the forward problem are naturally needed. Therefore, in this section, in particular, we present the necessary results of our previous work \cite{AS}, where the forward problem was studied.
In Section 4, we study an auxiliary problem related to the main problem. The proofs of the main results are contained in Section 5. The article ends with the section Conclusion.

\section{\textbf {Statement of the main results}}

Let $C[0, T]$ be the set of continuous functions defined on $[0,T]$ with the standard max-norm $||\cdot||_{C[0,T]}$ and let $C(H)=C([0,T]; H)$ stand
for a space of continuous $H$-valued functions $h(t)$ defined on $[0,T]$, and equipped with the norm
\[
||h||_{C(H)}=\max\limits_{0\leq t\leq T}||h(t)||.
\]

 When solving the inverse problem, we will investigate the Cauchy problem for various differential equations. In this case, by the solution of the problem we mean the classical solution, i.e. we will assume that all derivatives and functions involved in the equation are continuous with respect to the variable $t$. As an example, let us give the definition of the solution to the inverse problem.
\begin{defin}\label{def} A pair of functions $\{u(t), p(t)\}$ with the properties
$(D_{t}^{\rho })^{2}u(t)$, $Au(t)\in C((0,T]; H)$ and $D^{\rho}_{t}u(t),u(t)\in C(H)$, $p(t)\in C[0,T]$  and satisfying conditions
\eqref{prob1}, \eqref{ad}  is called \textbf{the
 solution} of the inverse problem.
\end{defin}

In order to formulate the main results of this paper, we introduce the Hilbert space of "smooth"
functions related to the degree of operator $A$.

Let $ \tau $ be an arbitrary real number. We introduce the power
of operator $ A $, acting in $ H $ according to the rule
\[
A^\tau h= \sum\limits_{k=1}^\infty \lambda_k^\tau h_k v_k, 
\]
here and everywhere below, by $h_k$  we will denote the Fourier coefficients of a vector $h\in H$: $h_k=(h,v_k)$.

Obviously, the domain of definition of this operator has the
form
\[
D(A^\tau)=\{h\in H:  \sum\limits_{k=1}^\infty \lambda_k^{2\tau}
|h_k|^2 < \infty\}.
\]
It is immediate from this definition that $D(A^\tau)\subseteq D(A^\sigma)$ for any $\tau\ge \sigma$.

On the set $D(A^\tau)$, we define the inner product
\[
(h,g)_\tau=\sum\limits_{k=1}^\infty \lambda_k^{2\tau} h_k \overline{g}_k =
(A^\tau h, A^\tau g)
\]

Then, $D(A^\tau)$ becomes a Hilbert
space with the norm $||h||^2_\tau=(h,h)_\tau.$
\begin{thm}\label{main2} Let  $\alpha>0$,  $\varphi_{0}\in H$, $\varphi_{1}\in D(A^{\frac{1}{2}})$ and $(D_t^{\rho})^{2} \psi(t)\in C[0,T]$. Further, let  $0<\epsilon<1$ be any fixed number and $q\in D(A^{1+\epsilon})$ and $f(t)\in C([0,T]; D(A^\epsilon))$, $Bq\neq 0$. Then the inverse problem has a unique solution $\{u(t), p(t)\}$.
\end{thm}
\begin{thm}\label{estimate} Let assumptions of Theorem~\ref{main2} be satisfied and let $\varphi_{0}\in D(A^{\frac{1}{2}}),\varphi_{1}\in D(A)$. Then  the solution to the invest problem obeys the stability estimate
\[\
   ||(D_{t}^{\rho})^{2} u||_{C(H)} +||D_t^\rho u||_{C(H)} + ||A u||_{C(H)} + ||p||_{C[0,T]}\leq C_{\rho, q, B, \epsilon} \bigg[  ||\varphi_{0}||_{\frac{1}{2}}+
  ||\varphi_{1}||_1 
  \]
  \begin{equation}\label{st}
  + ||\psi||_{C[0,T]} + +||D_t^\rho\psi||_{C[0,T]}+||(D_{t}^{\rho})^{2} \psi||_{C[0,T]}+ \max\limits_{0\leq t\leq T}||f(t)||_\epsilon\bigg],
 \end{equation}
where $C_{\rho, q, B, \epsilon}$ is a constant, depending only on $\rho, q$, $B$ and $\epsilon$.
\end{thm}

The choice of an abstract operator allows us to consider various known models. Note that the operator $A$ is only required to have a complete orthonormal system of eigenfunctions. Therefore, as $A$ one can consider any of the elliptic operators given in the work of Ruzhansky et al. \cite{Ruzhansky}. Let us consider the Laplace operator with the Dirichlet condition in more detail.

Let
$\Omega\subset \mathbb{R}^N$ be a bounded domain with sufficiently smooth
boundary $\partial \Omega$. Let $A_0$ denote the operator in
$L_2(\Omega)$ with domain $D(A_0)=\{g\in
C^2(\Omega)\cap C(\overline{\Omega}):\,g(x)=0,\,x\in\partial
\Omega\}$ and acting as $A_0g(x)=-\triangle g(x)$, where
$\triangle$ is the Laplace operator. Then (see, e.g.
\cite{Il}) $A_0$ has a complete in $L_2(\Omega)$ system of
orthonormal eigenfunctions $\{v_k(x)\}$ and a countable set of
nonnegative eigenvalues $\lambda_k$ ($\rightarrow +\infty$), and
$\lambda_1 = \lambda_1(\Omega)>0$.

Let $A$ stand for the operator, acting as $Ag(x) =\sum \lambda_k
g_k v_k(x)$ with the domain of definition $D(A)=\{ g\in
L_2(\Omega): \,\sum \lambda^2_k |g_k|^{2}<\infty\}$. Then it is not
hard to see, that $A$ is a positive self-adjoint extension in
$L_2(\Omega)$ of operator $A_0$. Accordingly, we can apply the Theorems
\ref{main2} and \ref{estimate} to  operator $A$ and hence to
inverse problem:
\begin{equation}\nonumber
\begin{cases}
  & (D_{t}^{\rho })^{2}u(x,t)+2\alpha D_{t}^{\rho }u(x,t)-\triangle u(x,t)=p(t)q(x)+f(x,t), \quad  0<t\le T; \\
& u(x, t)=0, \, x\in \partial \Omega, \quad 0<t\le T;\\
  
 & \underset{t\to 0}{\mathop{\lim }}\,D_{t}^{\rho }u(x,t)={{\varphi }_{0}}(x),\quad  x\in
 \Omega;     \\
 & u(x,0)={{\varphi }_{1}}(x), \quad x\in
\Omega,  \\
\end{cases}
\end{equation}
with an over-determination condition
\[
B[u(x,t)]\equiv \int\limits_{\Omega} u(x,t) dx=\psi(t), \quad
0\leq t\leq T.
\]
As a functional $B$ one can also consider $B[u(x,t)]\equiv u(x_0, t)$,
$x_0\in \overline{\Omega}$.

\section{\textbf {Preliminaries }}
\begin{defin}(Mittag-Leffler function) The function

\[
E_{\rho, \mu}(z)= \sum\limits_{k=0}^\infty \frac{z^k}{\Gamma(\rho
k+\mu)}
\]
is called the Mittag-Leffler function with two parametrs \cite{Gorenflo}, where $\rho>0$, $\mu \in \mathbb{C}$ and $z\in\mathbb{C}$.
\end{defin}
A further generalization of Mittag-Leffler function was introduced by Prabhakar \cite{Prah} as

\[
E^{\gamma}_{\rho, \mu}(z)= \sum\limits_{k=0}^\infty \frac{(\gamma)_{k}}{\Gamma(\rho k+\mu)}\cdot\frac{z^{k}}{k!},
\]
where $z \in \mathbb{C}$, $\rho$, $\mu$ and $\gamma$ are arbitrary positive constants, and $(\gamma)_{k}$ is the Pochhammer
symbol defined as 
\[
(\gamma)_{k}=\frac{\Gamma(\gamma+k)}{\Gamma(\gamma)}=\begin{cases}
  & 1,\quad  k=0, \gamma\neq0 ,\\
 & \gamma(\gamma+1)(\gamma+2)....(\gamma+k-1), \quad k\in\mathbb{N} .\\
\end{cases}
\] 
 \begin{lem}\label{ml}(asymptotic estimate \cite{Dzh66}, p. 133) Let  $\mu$  be an arbitrary complex number. Further let $\beta$ be a fixed number, such that $\frac{\pi}{2}\rho<\beta<\pi \rho$, and $\beta \leq |\arg z|\leq \pi$. Then the following asymptotic estimate holds
\[
E_{\rho, \mu}(z)= -\frac{z^{-1}}{\Gamma(\rho-\mu)} + O(|z|^{-2}),
\,\, |z|>1.
\]
\end{lem}
\begin{cor}\label{cor} Under the conditions of Lemma \ref{ml} one has
\[
|E_{\rho,\mu}(z)|\le \frac{C}{1+|z|}, \quad |z|\ge0,
\]
where  $C$-constant, independent of $z$.
\end{cor}

Everywhere below, by $C$ we will denote the constants not necessarily the same. And if a constant depends on, say, $A$, then we will denote it by $C_A$.

\begin{lem}\label{PDer}If $\rho>0$ and $\lambda \in \mathbb{C}$, then (see \cite{Pang})
\begin{equation}\nonumber
 D^{\rho}_{t}E_{\rho,1}(\lambda t^{\rho})=\lambda E_{\rho,1}(\lambda t^{\rho}) \quad t>0 ,
\end{equation}
\begin{equation}\nonumber
D^{\rho}_{t}\left(t^{\rho}E^{2}_{\rho,\rho+1}(\lambda t^{\rho}) \right)=E^{2}_{\rho,1}(\lambda t^{\rho}) \quad t>0 .
\end{equation}
\end{lem}
\begin{lem} The solution to the Cauchy problem
\begin{equation}\label{Dervitavi1}
\begin{cases}
  & D_{t}^{\rho }u(t)-\lambda u(t)=f(t),\quad  0<t\le T; \\
 & u(0)=0, \\
\end{cases}
\end{equation}
with $0<\rho<1$ and $\lambda\in\mathbb{C} $ has the form
\[
u(t)=\int_{0}^{t}(t-\tau)^{\rho-1}E_{\rho,\rho}(\lambda(t-\tau)^{\rho})f(\tau)d\tau.
\]
\end{lem}
The proof of this lemme for $\lambda\in\mathbb{R}$ can be found in~\cite{Kil} (p. 231). In a complex case, similar ideas will lead us to the same conclusion.

Consider the operator $E_{\rho, \mu} (t^{\rho} A): H\rightarrow H$ defined by the spectral theorem of J. von Neumann:
\[
E_{\rho, \mu} (t^{\rho} A)g(t) = \sum\limits_{k=1}^\infty E_{\rho,\mu} (t^{\rho}\lambda_{k}) g_k(t) v_k,
\]
here and everywhere below, by $g_k(t)$  we will denote the Fourier coefficients of a function $g(t)\in C(H)$: $g_k(t)=(g(t),v_k)$.

We will use the following Lemmas \ref{l}--\ref{forward} , proofs of which can be found in the paper \cite{AS}.

\begin{lem}\label{l} Let $g(t)\in C[0,T]$.Then the unique solution of the Cauchy problem
\begin{equation}\label{Integro}
\begin{cases}
  & D_{t}^{\rho }u(t)+2\lambda u(t)+\lambda^{2}J^{-\rho}_{t}u(t)=J^{-\rho}_{t}g(t), \quad  0<t\le T; \\
 & u(0)=0, \\
\end{cases}
\end{equation}
with $0<\rho<1$ and $\lambda\in \mathbb{C}$ has the form
\[
u(t)=\int_{0}^{t}(t-\tau)^{2\rho-1}E^{2}_{\rho,2\rho}(-\lambda(t-\tau)^{\rho})g(\tau)d\tau.
\]
\end{lem}

\begin{lem}\label{Estimate1} Let $\alpha>0$. Then for any $g(t)\in C(H)$  one has $E_{\rho,\mu}(-S t^{\rho})g(t)$, $SE_{\rho,\mu}(-S t^{\rho})g(t)\in C(H)$
and  $AE_{\rho,\mu}(-St^{\rho})g(t)\in C((0,T];H)$.  Moreover,the following estimate holds:
\begin{equation}\label{ES}
||E_{\rho, \mu} (-t^{\rho}S)g(t)||_{C(H)}\leq C||g(t)||_{C(H)}.
\end{equation}

If $g(t)\in D(A^{\frac{1}{2}})$ for all $t\in[0,T]$, then
\begin{equation}\label{SES1}
||SE_{\rho, \mu} (-t^{\rho} S)g(t)||_{C(H)} \leq C\max_{0\leq t \leq T}||g(t)||_{\frac{1}{2}},
\end{equation}
\begin{equation}\label{AES}
||AE_{\rho, \mu} (-t^{\rho} S)g(t)|| \leq Ct^{-\rho}\max_{0\leq t \leq T}||g(t)||_{\frac{1}{2}}, \quad t>0.
\end{equation}

Here $S$ has two states: $S^{-}$ and $S^{+}$,
\[
S^{-}=\alpha I-(\alpha^{2}I-A)^{\frac{1}{2}},\quad  S^{+}=\alpha I+(\alpha^{2}I-A)^{\frac{1}{2}}.
\]
\end{lem}
\begin{lem}\label{estimate2} Let   $\alpha>0$ and   $\lambda_{k}\neq \alpha^{2}$,  for all $k$. Then for  any $g(t)\in C(H)$  one has $R^{-1}E_{\rho,\mu}(-St^{\rho})g(t), SR^{-1}E_{\rho,\mu}(-St^{\rho})g(t)\in C(H)$ and $AR^{-1}E_{\rho, \mu} (-t^{\rho} S)g(t)\in C((0,T],H)$. Moreover,the following estimates hold:
\begin{equation}\label{RES}
||R^{-1}E_{\rho, \mu} (-t^{\rho} S)g(t)||_{C(H)} \leq C ||g(t)||_{C(H)},
\end{equation}
\begin{equation}\label{SRES}
||SR^{-1}E_{\rho, \mu} (-t^{\rho} S)g(t)||_{C(H)} \leq C ||g(t)||_{C(H)},
\end{equation}
\begin{equation}\label{ARES}
||AR^{-1}E_{\rho, \mu} (-t^{\rho} S)g(t)|| \leq Ct^{-\rho}||g(t)||_{C(H)}, \quad t>0.
\end{equation}

Here
\[
 R^{-1}=(\alpha^{2}I-A)^{-\frac{1}{2}}.
\]
\end{lem}
\begin{lem}\label{aep} Let $\alpha>0$ and $\lambda_{k}\neq \alpha^{2}$,  for all $k$. Then for any $g(t)\in C([0,T]; D(A^\epsilon))$, with  $0< \epsilon < 1 $ we have
\begin{equation}\label{A}
\bigg|\bigg|\int\limits_0^t(t-\tau)^{\rho-1} AR^{-1}E_{\rho, \rho}(-(t-\tau)^\rho S)g(\tau) d\tau   \bigg|\bigg|\leq C \max\limits_{0\leq t \leq T} ||g(t)||_\epsilon,
\end{equation}
\begin{equation}\label{SR}
\bigg|\bigg|\int\limits_0^t(t-\tau)^{\rho-1} SR^{-1}E_{\rho, \rho}(-(t-\tau)^\rho S)g(\tau) d\tau   \bigg|\bigg|\leq C \max\limits_{0\leq t \leq T} ||g(t)||_{\epsilon},
\end{equation}
\begin{equation}\label{R}
\bigg|\bigg|\int\limits_0^t(t-\tau)^{\rho-1}R^{-1}E_{\rho, \rho}(-(t-\tau)^\rho S)g(\tau) d\tau   \bigg|\bigg|\leq C \max\limits_{0\leq t \leq T} ||g(t)||_{\epsilon}.
\end{equation}
\end{lem}
\begin{lem}\label{Int} Let $\alpha>0$ and  $g(t)\in C(H)$. Then
  \[
  \bigg |\bigg|J^{-\rho}_{t}\bigg(\int_{0}^{t}(t-\tau)^{2\rho-1}E^{2}_{\rho,2\rho}(-\alpha(t-\tau)^{\rho})g(\tau)d\tau\bigg)\bigg|\bigg|_{C(H)}\leq C ||g(t)||_{C(H)},
\]
\begin{equation}\label{J1}
\bigg|\bigg|\int_{0}^{t}(t-\tau)^{2\rho-1}E^{2}_{\rho,2\rho}(-\alpha(t-\tau)^{\rho})g(\tau)d\tau\bigg|\bigg|_{C(H)}\leq C||g(t)||_{C(H)}.
\end{equation}
\end{lem}
\begin{lem}\label{forward}
   Let $\alpha>0$, $\varphi_{0}\in H$ and $\varphi_{1}\in D(A^{\frac{1}{2}})$.Further, let $\epsilon\in (0,1)$ be any fixed number and  $G(t)\in C([0,T]; D(A^{\epsilon}))$. Then the forward problem
   \begin{equation}\label{9}
\begin{cases}
  & (D_{t}^{\rho })^{2}u(t)+2\alpha D_{t}^{\rho }u(t)+Au(t)=G(t), \quad  0<t\le T; \\
 & \underset{t\to 0}{\mathop{\lim }}\,D_{t}^{\rho }u(t)={{\varphi }_{0}},\\
 & u(0)={{\varphi }_{1}}, \\
\end{cases}
\end{equation}
 has a unique solution.
\end{lem}
The solution to the problem \eqref{9} is considered in two cases:

Case I: $\forall k\in\mathbb{N}$ we have $\alpha^2\neq \lambda_k$;

Case II: $\exists k_0\in \mathbb{N}$, such that $\alpha^2=\lambda_{k_0}$.

With the same reasoning as in \cite{AS}, it is sufficient to work with the second case. For simplicity, we assume that there is only one $\lambda_{k_0}$ of this kind. Then the solution is (see \cite{AS})
\begin{equation}\label{f}\nonumber
u(t)=\frac{1}{2}\bigg[\tilde{E}_{\rho,1}(-S^{-}t^{\rho})+\tilde{E}_{\rho,1}(-S^{+}t^{\rho})\bigg]\varphi_{1}
+\frac{\alpha}{2}\bigg[R^{-1}\tilde {E}_{\rho,1}(-S^{-}t^{\rho})
\end{equation}
\[
-R^{-1}\tilde{E}_{\rho,1}(-S^{+}t^{\rho})\bigg]\varphi_{1}
+\frac{1}{2}\bigg[R^{-1}\tilde{E}_{\rho,1}(-S^{-}t^{\rho})-R^{-1}\tilde{E}_{\rho,1}(-S^{+}t^{\rho})\bigg]\varphi_{0}
\]
\[
+E_{\rho,1}(-\alpha t^{\rho})\varphi_{1k_0}v_{k_0}+\alpha t^{\rho} E^{2}_{\rho,1+\rho}(-\alpha t^{\rho})\varphi_{1k_0}v_{k_0}
\]
\[
+t^{\rho}E^{2}_{\rho,1+\rho}(-\alpha t^{\rho})\varphi_{0k}v_{k_0}+\int_{0}^{t}(t-\tau)^{2\rho-1}E^{2}_{\rho,2\rho}(-\alpha(t-\tau)^{\rho})G_{k_0}(\tau)d\tau v_{k_0}
\]
\[
+\frac{1}{2}\int_{0}^{t}(t-\tau)^{\rho-1}\bigg[R^{-1}\tilde{E}_{\rho,\rho}(-S^{-}(t-\tau)^{\rho})-R^{-1}\tilde{E}_{\rho,\rho}(-S^{+}(t-\tau)^{\rho})\bigg]G(\tau)d\tau,
\]
where we denote  $$\varphi_{jk}=(\varphi_{j},v_{k}) $$
and
$$\tilde{E}_{\rho,\mu}(-St^\rho)g=\sum\limits_{k\neq k_0} E_{\rho,\mu} (-(\alpha\pm\sqrt{\alpha^{2}-\lambda_{k}})t^{\rho}) g_k v_k.$$

\begin{lem}\label{3.13}Let $\alpha>0$. Then for any $g(t)\in D(A)$  one has 
 $AE_{\rho,\mu}(-St^{\rho})g(t)\in C(H)$.  Moreover, the following estimate holds:
\begin{equation}\label{AES1}
||A E_{\rho, \mu} (-t^{\rho} S)g(t)||_{C(H)}\leq C\max_{0\leq t \leq T}||g(t)||_{1}.    
\end{equation}

Here $S$ has two states: $S^{-}$ and $S^{+}$,
\[
S^{-}=\alpha I-(\alpha^{2}I-A)^{\frac{1}{2}},\quad  S^{+}=\alpha I+(\alpha^{2}I-A)^{\frac{1}{2}}.
\]

\end{lem}
\proof
By using Parseval's equality, one has
\[
||AE_{\rho, \mu} (-S^{-} t^{\rho})g(t)||^2=\sum\limits_{k=1}^\infty \left|\lambda_{k}E_{\rho, \mu} \left(-\left(\alpha-\sqrt{\alpha^{2}-\lambda_{k}}\right)t^{\rho}\right) g_k(t)\right|^2.
\]
According to Corollary \ref{cor}, finally we have
\[
||AE_{\rho, \mu} (-S^{-}t^{\rho})g(t)||^2\leq C \sum\limits_{k=1}^\infty \bigg|\frac{ \lambda_{k}g_k(t)}{1+|\alpha-\sqrt{\alpha^{2}-\lambda_{k}}|t^{\rho}}\bigg|^2\leq C ||g(t)||_{1}^2.
 \]
Similar estimate is proven in precisely the same way for  operator $S^{+}$.

\endproof
\begin{lem}
Let   $\alpha>0$ and   $\lambda_{k}\neq \alpha^{2}$,  for all $k$. Then for  any $g(t)\in D(A^{\frac{1}{2}})$  one has $AR^{-1}E_{\rho, \mu} (-t^{\rho} S)g(t)\in C(H)$. Moreover, the following estimate holds:

\begin{equation}\label{ARES1}
 ||AR^{-1}E_{\rho, \mu} (-t^{\rho} S)g(t)||_{C(H)} \leq C  \max_{0\leq t \leq T}||g(t)||_{\frac{1}{2}},   
\end{equation}
\end{lem}
\proof
In proving the lemma, we use Parseval's equality and Corollary \ref{cor} similarly to the proof of  Lemma \ref{3.13}:
\[
||AR^{-1}E_{\rho, \mu} (-t^{\rho}S^{-})g(t)||^{2}\leq C \sum\limits_{k=1}^\infty \bigg|\frac{1}{\sqrt{\alpha^{2}-\lambda_{k}}}\frac{\lambda_{k}g_k(t)}{1+t^{\rho}|\alpha-\sqrt{\alpha^{2}-\lambda_{k}}|}\bigg|^{2},
\]
\[
u_{\lambda_{k}}(t)=\frac{\lambda_{k}^{2}|g_k(t)|^{2}}{|\sqrt{\alpha^{2}-\lambda_{k}}|^{2}|(1+t^{\rho}|\alpha-\sqrt{\alpha^{2}-\lambda_{k}|})^{2}} \quad  \underset{\lambda_{k}\to \infty}{\sim} \quad  v_{\lambda_{k}}(t)=\frac{\lambda_{k}^{2}|g_k(t)|^{2}}{\lambda_{k}(1+t^{\rho}\sqrt{\lambda_{k}})^{2}},
\]
\[
\sum_{k=1}^{\infty}v_{\lambda_{k}}(t)=\sum_{k=1}^{\infty}\frac{ \lambda_{k}|g_k(t)|^{2}}{(1+t^{\rho}\sqrt{\lambda_{k}})^{2}}\leq ||g(t)||^{2}_{\frac{1}2}.
\]
We have used the notation here:  ${u_{\lambda_{k}}} \underset{\lambda_{k}\to \infty}{\sim} {v_{\lambda_{k}}}$ means 
$
\lim\limits_{\lambda_{k}\to \infty}\frac{{u_{\lambda_{k}}}}{{v_{\lambda_{k}}}}=1 $.

Similar estimate is proven in precisely the same way for  operator $S^{+}$.
\endproof

\section{Auxiliary Problem}
Let us consider the following auxiliary Cauchy problem
\begin{equation}\label{Cauchy3}
\begin{cases}
  (D_{t}^{\rho })^{2}w(t)+2\alpha D_{t}^{\rho }w(t)+Aw(t)=-\mu(t)Aq-2\alpha qD^{\rho}_{t}\mu(t)+f(t), \quad 0< t\le T; \\
 \underset{t\to 0}{\mathop{\lim }}\,D_{t}^{\rho }w(t)={{\varphi }_{0}},\\
 w(0)={{\varphi }_{1}}, \\
\end{cases}
\end{equation}
where the scalar function $\mu(t)$ is the unique solution of the Cauchy problem (see, e.g., \cite{AS}):
\begin{equation}\label{mu}
\begin{cases}
  (D_{t}^{\rho })^{2}\mu(t)=p(t), \quad  0< t\le T; \\
 \underset{t\to 0}{\mathop{\lim }}\,D_{t}^{\rho }\mu(t)=0,\\
 \mu(0)=0, \\
\end{cases}
\end{equation}

\begin{lem}\label{lem} Let $w(t)$ be a solution of problem \eqref{Cauchy3} and $Bq\neq 0$. Then the unique solution $\{u(t),p(t)\}$ to the inverse problem  \eqref{prob1},\eqref{ad} has the form
 \begin{equation}\label{wq}
   u(t)=w(t)+\mu(t)q,
 \end{equation}
 \begin{equation}\label{p}
   p(t)=\frac{1}{Bq}\bigg[(D^{\rho}_{t})^{2}\psi(t)-B[(D^{\rho}_{t})^{2}w(t)]\bigg],
 \end{equation}
 where
 \begin{equation}\label{mub}
   \mu(t)=\frac{1}{Bq}\bigg[\psi(t)-B[w(t)]\bigg].
 \end{equation}
\end{lem}
\proof Substitute the function  $u(t)$, defined by equality \eqref{wq} into equation \eqref{prob1}. Then
\[
(D^{\rho}_{t})^{2}w(t)+(D^{\rho}_{t})^{2}\mu(t)q+2\alpha D^{\rho}_{t}w(t)+2\alpha D^{\rho}_{t}\mu(t)q+Aw(t)+\mu(t)Aq=p(t)q+f(t).
\]
Since $(D^{\rho}_{t})^{2}\mu(t)=p(t)$ we obtain equation \eqref{Cauchy3}. As for the initial conditions, again by virtue of \eqref{mu} we get
\[
\lim_{t\to 0}D^{\rho}_{t}u(t)=\lim_{t\to 0}D^{\rho}_{t}w(t)+\lim_{t\to 0}D^{\rho}_{t}\mu(t)q=\varphi_{0},
\]
\[
u(0)=w(0)+\mu(0)q=w(0)=\varphi_{1}.
\]
Thus function $u(t)$, defined as \eqref{wq} is a solution of the Cauchy problem \eqref{prob1}.

Let us prove equation \eqref{p}. Apply \eqref{wq} to obtain
\[
\mu(t)q=w(t)-u(t).
\]
On the other hand, due to equality \eqref{ad}, we have
   \begin{equation}\nonumber
     \mu(t) Bq =B[w(t)] - \psi(t),
     \end{equation}
    or, since $Bq\neq 0$, we get \eqref{mub}). Finally, using equality $(D_t^{\rho})^{2} \mu(t)= p(t)$, we obtain
    \begin{equation}\nonumber
    p(t)=\frac{1}{B q}\bigg[(D_t^{\rho})^{2} \psi(t)- B[(D_t^{\rho})^{2} w(t)]\bigg],
\end{equation}
    which coincides with (\ref{p}).
\endproof

Thus, to solve the inverse problem \eqref{prob1}, \eqref{ad}, it is sufficient to solve the Cauchy problem \eqref{Cauchy3}.

\begin{thm}\label{u} Under the assumptions of Theorem \ref{main2}, the problem \eqref{Cauchy3}, with $\mu$ defined in (\ref{mub}), has a unique solution.
\end{thm}
\proof
 Suppose that 
\[
    G(t)=-\frac{Aq}{Bq} (\psi(t)-B[w(t)])-\frac{2\alpha q}{Bq}(D^{\rho}_{t}\psi(t)-B[D^{\rho}_{t}w(t)]) 
   \]
\begin{equation}\label{Gsep}
    +f(t) \in C([0,T]; D(A^\epsilon))
\end{equation}
with some $\epsilon \in (0,1)$. Then by Lemma \ref{forward} the solution $w(t)$ to the problem \eqref{Cauchy3} solves the following integral equation
   \begin{equation}\label{f}
w(t)=\frac{1}{2}\bigg[\tilde{E}_{\rho,1}(-S^{-}t^{\rho})+\tilde{E}_{\rho,1}(-S^{+}t^{\rho})\bigg]\varphi_{1}+\frac{\alpha}{2}\bigg[R^{-1}\tilde {E}_{\rho,1}(-S^{-}t^{\rho})
\end{equation}
\[
-R^{-1}\tilde{E}_{\rho,1}(-S^{+}t^{\rho})\bigg]\varphi_{1}
+\frac{1}{2}\bigg[R^{-1}\tilde{E}_{\rho,1}(-S^{-}t^{\rho})-R^{-1}\tilde{E}_{\rho,1}(-S^{+}t^{\rho})\bigg]\varphi_{0}
\]
\[
+E_{\rho,1}(-\alpha t^{\rho})\varphi_{1k_0}v_{k_0}+\alpha t^{\rho} E^{2}_{\rho,1+\rho}(-\alpha t^{\rho})\varphi_{1k_0}v_{k_0}
\]
\[
+t^{\rho}E^{2}_{\rho,1+\rho}(-\alpha t^{\rho})\varphi_{0k}v_{k_0}+\int_{0}^{t}(t-\tau)^{2\rho-1}E^{2}_{\rho,2\rho}(-\alpha(t-\tau)^{\rho})G_{k_0}(\tau)v_{k_0}d\tau
\]
\[
+\frac{1}{2}\int_{0}^{t}(t-\tau)^{\rho-1}\bigg[R^{-1}\tilde{E}_{\rho,\rho}(-S^{-}(t-\tau)^{\rho})-R^{-1}\tilde{E}_{\rho,\rho}(-S^{+}(t-\tau)^{\rho})\bigg]G(\tau)d\tau.
\]

 Rewrite this equation as
\begin{equation}\label{VE2}
        w(t)=F(t)+\frac{1}{2}\int_{0}^{t}(t-\tau)^{\rho-1}R^{-1}\tilde{E}_{\rho,\rho}(-S^{-}(t-\tau)^{\rho})B_{1}[w(\tau)]d\tau
     \end{equation}
     \[
        -\frac{1}{2}\int_{0}^{t}(t-\tau)^{\rho-1}R^{-1}\tilde{E}_{\rho,\rho}(-S^{+}(t-\tau)^{\rho}B_{1}[w(\tau)]d\tau
        \]
        \[
        +\int_{0}^{t}(t-\tau)^{2\rho-1}E^{2}_{\rho,2\rho}(-\alpha(t-\tau)^{\rho})B_{1}[w_{k_0}(\tau)]d\tau v_{k_0},
     \]
 where
   \[
F(t)=\frac{1}{2}\big[\tilde{E}_{\rho,1}(-S^{-}t^{\rho})+\tilde{E}_{\rho,1}(-S^{+}t^{\rho})\big]\varphi_{1}+\frac{\alpha}{2}R^{-1}\big[\tilde{E}_{\rho,1}(-S^{-}t^{\rho})-\tilde{E}_{\rho,1}(-S^{+}t^{\rho})\big]\varphi_{1}
\]
\[
+\frac{1}{2}R^{-1}\big[\tilde{E}_{\rho,1}(-S^{-}t^{\rho})-\tilde{E}_{\rho,1}(-S^{+}t^{\rho})\big]\varphi_{0}
\]
\[
+E_{\rho,1}(-\alpha t^{\rho})\varphi_{1k_{0}}v_{k_{0}}+\alpha t^{\rho}E^{2}_{\rho,1+\rho}(-\alpha t^{\rho})\varphi_{1k_{0}}v_{k_{0}}+t^{\rho}E^{2}_{\rho,1+\rho}(-\alpha t^{\rho})\varphi_{0k_{0}}v_{k_{0}}
\]
\[
+\int_{0}^{t}(t-\tau)^{2\rho-1}E^{2}_{\rho,2\rho}(-\alpha(t-\tau)^{\rho})(-\frac{\lambda_{k_{0}}q_{k_{0}}}{Bq
}\psi(\tau)-\frac{2\alpha q_{k_{0}}}{Bq}D^{\rho}_{\tau}\psi(\tau)+f_{k_{0}}(\tau))v_{k_{0}}d\tau
\]
\[
+\frac{1}{2}\int_{0}^{t}(t-\tau)^{\rho-1}R^{-1}\tilde{E}_{\rho,\rho}(-S^{-}(t-\tau)^{\rho})(-\frac{Aq}{Bq}\psi(\tau)-\frac{2\alpha q}{Bq}D^{\rho}_{\tau}\psi(\tau)+f(\tau))d\tau
\]
\[
-\frac{1}{2}\int_{0}^{t}(t-\tau)^{\rho-1}R^{-1}\tilde{E}_{\rho,\rho}(-S^{+}(t-\tau)^{\rho})(-\frac{Aq}{Bq}\psi(\tau)-\frac{2\alpha q}{Bq}D^{\rho}_{\tau}\psi(\tau)+f(\tau))d\tau,
\]
  and
    \[
    B_1 [w(\tau)]=\frac{Aq}{Bq}B [w(\tau)]+\frac{2\alpha q}{Bq}B [D^{\rho}_{\tau}w(\tau)].
    \] 
    
Here our aim is to show the following two facts:

    (a) There exists a unique solution to the integral equation \eqref{VE2};
    
    (b) The solution to \eqref{VE2} is a solution to  equation \eqref{Cauchy3}.

To achieve these goals, we will use the original ideas of work Ashyralyev and Al-Hazaimeh \cite{Ashyr3}.
    
We can rewrite the derivative $D^{\rho}_{t}w(t)$  using the equations \eqref{Dervitavi1}, \eqref{Integro}:
\[
\frac{2\alpha q}{Bq} D^{\rho}_{t}w(t)=\frac{2\alpha q}{Bq}D^{\rho}_{t}F(t)-\frac{\alpha q}{Bq}\int_{0}^{t}(t-\tau)^{\rho-1}R^{-1}S^{-}\tilde{E}_{\rho,\rho}(-S^{-}(t-\tau)^{\rho})B_{1}[w(\tau)]d\tau
\]
\[
+\frac{\alpha q}{Bq}\int_{0}^{t}(t-\tau)^{\rho-1}R^{-1}S^{+}\tilde{E}_{\rho,\rho}(-S^{+}(t-\tau)^{\rho})B_{1}[w(\tau)]d\tau
\]
\[
-\frac{2\alpha^{2} q}{Bq}\int_{0}^{t}(t-\tau)^{2\rho-1}E^{2}_{\rho,2\rho}(-\alpha(t-\tau)^{\rho})B_{1}[w_{k_0}(\tau)]d\tau v_{k_0}
\]
\[
-\frac{2\alpha^{3} q}{Bq}J^{-\rho}_{t}\left[\int_{0}^{t}(t-\tau)^{2\rho-1}E^{2}_{\rho,2\rho}(-\alpha(t-\tau)^{\rho})B_{1}[w_{k_0}(\tau)]d\tau v_{k_0}\right]+\frac{2\alpha q}{Bq}J^{-\rho}_{t}B_{1}[w_{k_{0}}(t)]v_{k_{0}}.
\]
 We use Lemma \ref{PDer} and \eqref{Integro}, \eqref{Dervitavi1}  to calculate $D^{\rho}_{t}F(t)$,
 \[
D^{\rho}_{t}F(t)=\frac{1}{2}\left[-S^{-}\tilde{E}_{\rho,1}(-S^{-}t^{\rho})-S^{+}\tilde{E}_{\rho,1}(-S^{+}t^{\rho})\right]\varphi_{1}-\frac{\alpha}{2}R^{-1}S^{-}\tilde{E}_{\rho,1}(-S^{-}t^{\rho})\varphi_{1}
\]
\[
+\frac{\alpha}{2}R^{-1}S^{+}\tilde{E}_{\rho,1}(-S^{+}t^{\rho})\varphi_{1}+\frac{1}{2}\left[-R^{-1}S^{-}\tilde{E}_{\rho,1}(-S^{-}t^{\rho})+R^{-1}S^{+}\tilde{E}_{\rho,1}(-S^{+}t^{\rho})\right]\varphi_{0}
\]
\[
-\alpha E_{\rho,1}(-\alpha t^{\rho})\varphi_{1k_{0}}v_{k_{0}}+\alpha E^{2}_{\rho,1+\rho}(-\alpha t^{\rho})\varphi_{1k_{0}}v_{k_{0}}+E^{2}_{\rho,1+\rho}(-\alpha t^{\rho})\varphi_{0k_{0}}v_{k_{0}}
\]
\[
-2\alpha\int_{0}^{t}(t-\tau)^{2\rho-1}E^{2}_{\rho,2\rho}(-\alpha(t-\tau)^{\rho})(-\frac{\lambda_{k_{0}}q_{k_{0}}}{Bq}\psi(\tau)-\frac{2\alpha q_{k_{0}}}{Bq}D^{\rho}_{\tau}\psi(\tau)+f_{k_{0}}(\tau))v_{k_{0}}d\tau
\]
\[
-\alpha^{2}J^{-\rho}_{t}\left(\int_{0}^{t}(t-\tau)^{2\rho-1}E^{2}_{\rho,2\rho}(-\alpha(t-\tau)^{\rho})(-\frac{\lambda_{k_{0}}q_{k_{0}}}{Bq}\psi(\tau)-\frac{2\alpha q_{k_{0}}}{Bq}D^{\rho}_{\tau}\psi(\tau)+f_{k_{0}}(\tau))v_{k_{0}}d\tau\right)
\]
\[
+J^{-\rho}_{t}\left(-\frac{\lambda_{k_{0}}q_{k_{0}}}{Bq}\psi(t)-\frac{2\alpha q_{k_{0}}}{Bq}D^{\rho}_{t}\psi(t)+f_{k_{0}}(t)\right)v_{k_{0}}
\]
\[
-\frac{1}{2}\int_{0}^{t}(t-\tau)^{\rho-1}R^{-1}S^{-}\tilde{E}_{\rho,\rho}(-S^{-}(t-\tau)^{\rho})(-\frac{Aq}{Bq}\psi(\tau)-\frac{2\alpha q}{Bq}D^{\rho}_{\tau}\psi(\tau)+f(\tau))d\tau
\]
\[
+\frac{1}{2}\int_{0}^{t}(t-\tau)^{\rho-1}R^{-1}S^{+}\tilde{E}_{\rho,\rho}(-S^{+}(t-\tau)^{\rho})(-\frac{Aq}{Bq}\psi(\tau)-\frac{2\alpha q}{Bq}D^{\rho}_{\tau}\psi(\tau)+f(\tau))d\tau,
\]
and we derive the following integral equation from  equation \eqref{VE2}
\begin{equation}\label{Volterra}
z(t)=\Phi(t)+\int_{0}^{t}K(t,\tau)B[z(\tau)]d\tau,
\end{equation}
where

\begin{equation}\label{zw}
    z(t)=\frac{Aq}{Bq}w(t)+\frac{2\alpha q}{Bq}D^{\rho}_{t}w(t),
\end{equation}

\[
\Phi(t)=\frac{Aq}{Bq}F(t)+\frac{2\alpha q}{Bq}D^{\rho}_{t}F(t)
\]
and
\[
K(t,\tau)=\frac{Aq}{2Bq}(t-\tau)^{\rho-1}R^{-1}\left
(E_{\rho,\rho}(-S^{-}(t-\tau)^{\rho})-E_{\rho,\rho}(-S^{+}(t-\tau)^{\rho}) \right)
\]
\[
+\frac{Aq}{Bq}(t-\tau)^{2\rho-1}E^{2}_{\rho,2\rho}(-\alpha(t-\tau)^{\rho})-\frac{\alpha q}{Bq}(t-\tau)^{\rho-1}R^{-1}S^{-}E_{\rho,\rho}(-S^{-}(t-\tau)^{\rho})
\]
\[
+\frac{\alpha q}{Bq}(t-\tau)^{\rho-1}R^{-1}S^{+}E_{\rho,\rho}(-S^{+}(t-\tau)^{\rho})+\frac{2\alpha q}{Bq \Gamma(\rho)}(t-\tau)^{\rho-1}
\]
\[
-\frac{2\alpha^{3} q}{Bq \Gamma(\rho)}(t-\tau)^{\rho-1}\int_{0}^{\tau}(\tau-\xi)^{2\rho-1}E^{2}_{\rho,2\rho}(-\alpha(\tau-\xi)^{\rho})d\xi.
\]
From estimates \eqref{RES}, \eqref{SRES} and Lemma \ref{Int}, we obtain the following estimate for arbitrary $g(t)\in C(H)$
\begin{equation}\label{kernelestimate}
\left|\left|\int_{\xi}^{\eta} K(t,\tau)g(\tau)d\tau \right|\right|\leq C_{B,q,\rho}\frac{(\eta-\xi)^{\rho}}{\rho}||g(t)||_{C([\xi,\eta],H)},
\end{equation}
where $0\le\xi<\eta\le T$.
\begin{rem}\label{fcontinuous}
    Notice that, Lemmas \ref{Estimate1}--\ref{Int} imply that $F(t)$ and $D_t^\rho F(t)$ are of class $C(H)$.
\end{rem}

Consider the Volterra integral equation \eqref{Volterra}. According to Remark \ref{fcontinuous}, we obtain $\Phi(t)\in C(H)$. In this case we have the following lemma.

\begin{lem}\label{z}
There is a unique solution $z(t)\in C(H)$  to the integral equation (\ref{Volterra}).
\end{lem}
\proof \textit{of the Lemma \ref{z}}.
Since $B: H\rightarrow \mathbb{R}$ is a bounded linear functional, let $b$ be its norm: $|B \varphi(t)|\leq b ||\varphi||_{C(H)}$, for all $\varphi\in C(H)$. Equation \eqref{Volterra} makes sense in any interval $[0, t_1]\in [0,T]$, $(0<t_1<T)$. One can find $t_1$ such that
\begin{equation}\label{t1}
C_{B,q,\rho} b \frac{t_1^\rho}{\rho}<1
\end{equation}
and show the existence of a unique solution $z(t)\in C([0, t_1]; H)$ to  equation \eqref{Volterra} on the interval $[0, t_1]$. For this we use the Banach fixed point theorem for the space $C([0, t_1]; H)$ (see, for example, \cite{Kil}, theorem 1.9, p. 68), with the distance given by
\[
d(z_1, z_2) = \max\limits_{t\in [0,t_1]} ||z_1(t)- z_2(t)||.
\]
We denote the right-hand  side of the equation \eqref{Volterra} by $Pz(t)$, where $P$ is the corresponding linear operator. To apply the Banach fixed point theorem we have to prove the following:

(1) if $z(t)\in C([0, t_1]; H)$, then $Pz(t)\in C([0, t_1]; H)$;

(2) for any $z_1, z_2 \in C([0, t_1]; H)$ one has
\[
\max\limits_{t\in [0,t_1]}||Pz_1-Pz_2||\leq \delta \max\limits_{t\in [0,t_1]} ||z_1-z_2||, \,\, \delta<1.
\]

 Remark \ref{fcontinuous} and \ref{kernelestimate} imply condition (1).
 
 Note that \eqref{kernelestimate} is equivalent to
\[
\left|\left|\int_{0}^{t_{1}}K(t,\tau)B[z_{1}(\tau)-z_{2}(\tau)]d\tau\right|\right|\leq \delta \max_{0\leq t\leq t_{1}}||z_{1}(t)-z_{2}(t)||,
\]
 where $\delta=C_{B,q,\rho}b\frac{t_{1}^{\rho}}{\rho}<1$ since condition \eqref{t1}.

Hence by the Banach fixed point theorem, there exists a unique solution $z^\star (t)\in C([0, t_1]; H)$ to equation \eqref{Volterra} on the interval $[0, t_1]$ and this solution is a limit of the convergent
sequence $z_n(t)=P^n \Phi(t)= P P^{n-1} \Phi(t)$:
\[
\lim\limits_{n\rightarrow \infty} \big[ \max\limits_{t\in [0, t_1]} ||z_n(t)-z^\star(t)||\big]=0.
\]
Since  estimate \eqref{kernelestimate} we can continue this argument beyond $t_1$. Next we consider the interval $[t_1, 2t_{1}]$. Rewrite  equation \eqref{Volterra} in the form
\[
z(t)=\Phi_1(t)+\int\limits_{t_1}^t K(t,\tau)B[z(\tau)]d\tau,
\]
where
\[
\Phi_1(t)=\Phi(t)+ \int\limits_{0}^{t_1}K(t,\tau)[Bz(\tau)]d\tau,
\]
is a known function, since the function $z(t)$ is uniquely defined on the interval $[0, t_1]$. Using the same arguments as above, we derive that there exists a unique solution $z^\star(t)\in C([t_1, 2t_1]; H)$ to equation \eqref{Volterra} on the interval $[t_1, 2t_1].$ It is also easy to see that $z(t_1)=z^*(t_1)$, so $z^*(t)$ is a continuation of $z(t)$ to the interval $[t_1, 2t_1]$. 
Repeating this process several times, we conclude that there exists a unique solution $z^\star(t)\in C([0, T]; H)$ to equation \eqref{Volterra} on the interval $[0, T]$, and this solution is a limit of the convergent
sequence $z_n(t)\in C([0,T]; H)$:
\begin{equation}\label{wn}\nonumber
\lim\limits_{n\rightarrow \infty} ||z_n(t)-z^\star(t)||_{C(H)}=0,
\end{equation}
with the choice of certain $z_n$ on each $[0, t_1], \cdots [nt_1, T]$.
\endproof

\textit{Continuation of the proof of Theorem \ref{u}}. 
Using the obtained solution $z(t)\in C(H)$ we go back to  equation ( \eqref{zw}, considering it as a Cauchy problem with initial data $w(0)=\varphi_1$. By the existence and uniqueness result (see \cite{Kil}, theorem 4.3, p. 231) we have a unique solution $w(t)$ of the considering Cauchy problem. This concludes the existence of the solution to the integral equation \eqref{VE2}.To prove uniqueness we use contradiction argument. Suppose there are two solutions $w_1(t)$ and $w_2(t)$ to the integral equation \eqref{VE2}. Then we know the corresponding functions $z_1(t)$ and $z_2(t)$ defined as \eqref{zw}. Thus both of them are solutions to the Volterra equation \eqref{Volterra} and by the uniqueness shown above we have that $z_1(t)\equiv z_2(t)$. Using the initial values we get $w_1(t)\equiv w_2(t)$. Thus we finish the part (a) of our aim.

Now we need to show that the solution $w(t)$ to  problem \eqref{VE2} satisfies equation \eqref{Cauchy3}. In order to show that we should have $(D_t^\rho)^{2} w(t),\,D_t^\rho w(t),\, Aw(t)\in C((0, T]; H)$. We already know $ D_t^\rho w(t)\in C((0, T]; H)$ by \eqref{zw}.

We start with $Aw(t)$ and for this purpose we consider $AF(t)$ and note that, by estimates \eqref{AES}, \eqref{ARES}, Corollary \ref{cor}, \eqref{A}, \eqref{J1} we obtain
\begin{equation}\label{aft}
||AF(t)||\leq C( t^{-\rho}||\varphi_{1}||_{\frac{1}{2}}+t^{-\rho}||\varphi_{0}||
+1)+\frac{|q_{k_{0}}|\alpha^{2}}{||Bq||}||D^{\rho}_{t}\psi(t)||_{C[0,T]}
\end{equation}
\[
+C_{\rho}\bigg[\frac{|q_{k_{0}}|\alpha^{2}}{||Bq||}||\psi(t)||_{C[0,T]}+C\max_{0\leq t \leq T}|f_{k_{0}}(t)|
\]
\[
+\frac{||A^{1+\epsilon}q||}{||Bq||}||\psi(t)||_{C[0,T]}
+\frac{2\alpha||q||}{||Bq||}||D^{\rho}_{t}\psi(t)||_{C[0,T]}]+C\max_{0\leq t \leq T}||f(t)||_{\epsilon}\bigg],
\]
for every $t\in(0,T]$.

Moreover we have the following estimate for the integral part of $Aw(t)$:
\begin{equation}\label{awt1}
\bigg|\bigg|\frac{1}{2}\int_{0}^{t}(t-\tau)^{\rho-1}AR^{-1}\tilde{E}_{\rho,\rho}(-S^{-}(t-\tau)^{\rho})B_{1}[w(\tau)]d\tau
\end{equation}
\[
-\frac{1}{2}\int_{0}^{t}(t-\tau)^{\rho-1}AR^{-1}\tilde{E}_{\rho,\rho}(-S^{+}(t-\tau)^{\rho})B_{1}[w(\tau)]d\tau
\]
\[
+\int_{0}^{t}(t-\tau)^{2\rho-1}\alpha^{2}E^{2}_{\rho,2\rho}(-\alpha(t-\tau)^{\rho})B_{1}[w_{k_0}(\tau)]d\tau v_{k_{0}}\bigg|\bigg|\leq C \max_{0\leq t\leq T}\bigg|\bigg|\frac{Aq}{Bq}B[w(t)]\bigg|\bigg|_{\epsilon}
\]
\[
+C \max_{0\leq t\leq T}\bigg|\bigg|\frac{2\alpha q}{Bq}B[D^{\rho}_{t}w(t)]\bigg|\bigg|_{\epsilon}+C_{\rho}\max_{0\leq t\leq T}\bigg[\frac{\alpha^{2}|q_{k_{0}}|}{||Bq||}|B[w_{k_{0}}(t)]|
+\frac{2\alpha |q_{k_{0}}|}{||Bq||}\left
|B[D^{\rho}_{t}w_{k_{0}}(t)]\right|\bigg]
\]
\[
\leq C_{B,q}\bigg[||A^{1+\epsilon}q||||w(t)||_{C(H)}+2\alpha||q||_{\epsilon}||D^{\rho}_{t}w(t)||_{C(H)}
\]
\[
+\alpha^{2}|q_{k_{0}}|||w_{k_{0}}(t)||_{C[0,T]}+2\alpha |q_{k_{0}}|||D^{\rho}_{t}w_{k_{0}}(t)||_{C[0,T]}  \bigg].
\]
Thus, $Aw(t)\in C((0, T]; H)$.

We have that  $w(t)$ and $D_t^\rho w(t)$ are of class $C(H)$. 
We set 
\[
    G(t)=-\frac{Aq}{Bq} (\psi(t)-B[w(t)])-\frac{2\alpha q}{Bq}(D^{\rho}_{t}\psi(t)-B[D^{\rho}_{t}w(t)])+f(t) 
   \]
By assumptions of the Theorem \ref{main2} we know $G(t)\in C([0,T], D(A^\epsilon)) $.Therefore, the  Cauchy problem 
\begin{equation*}
    \begin{cases}
  & (D_{t}^{\rho })^{2}\tilde{w}(t)+2\alpha D_{t}^{\rho }\tilde{w}(t)+A\tilde{w}(t)=G(t), \quad 0< t\le T; \\
 & \underset{t\to 0}{\mathop{\lim }}\,D_{t}^{\rho }\tilde{w}(t)={{\varphi }_{0}},\\
 & \tilde{w}(0)={{\varphi }_{1}}, \\
\end{cases}
\end{equation*}
has a solution given by the expression the same as in  equation \eqref{f}. By the definition of $G(t)$ the right-hand side of \eqref{f} is exactly right-hand side of \eqref{VE2} in terms of $\tilde{w}(t)$ and our $w(t)$ is the solution to \eqref{VE2}. Therefore $\tilde{w}(t)=w(t)$. 
\endproof
\begin{rem}\label{rem1}
Using estimates (\ref{AES1}), (\ref{ARES1}) and similar ideas as in the proof of  estimate \eqref{aft} we have the following estimate in the closed interval $[0,T]$
\[
 ||AF(t)||\leq C\big(||\varphi_{1}||_{1}+||\varphi_{0}||_{\frac{1}{2}}\big)+
C_{\rho}\bigg[\frac{|q_{k_{0}}|\alpha^{2}}{||Bq||}||\psi(t)||_{C[0,T]}+\frac{|q_{k_{0}}|\alpha^{2}}{||Bq||}||D^{\rho}_{t}\psi(t)||_{C[0,T]}
\]
\[
+C\max_{0\leq t \leq T}|f_{k_{0}}(t)|+\frac{||A^{1+\epsilon}q||}{||Bq||}||\psi(t)||_{C[0,T]}+\frac{2\alpha||q||}{||Bq||}||D^{\rho}_{t}\psi(t)||_{C[0,T]}
\]
\[
+C\max_{0\leq t \leq T}||f(t)||_{\epsilon}\bigg].
\]
\end{rem}

\section{Proof of Theorem \ref{estimate}}
First we prove the following statement on the stability of the solution to problem \eqref{Cauchy3}, \eqref{mu}.

\begin{thm}\label{estimate3} Let assumptions of Theorem \ref{estimate} be satisfied. Then  the solution to problem (\ref{Cauchy3}), (\ref{mu}) obeys the stability estimate
    \begin{equation}\label{se}
        ||(D_t^\rho)^{2} w(t)||_{C(H)}\leq C_{\rho, q, B, \epsilon} \big[ ||\varphi_{0}||_\frac{1}{2}+||\varphi_{1}||_{1}+ ||\psi(t)||_{C[0,T]}
       \end{equation}
  \[
        +||D^{\rho}_{t}\psi(t)||_{C[0,T]}+ \max\limits_{0\leq t\leq T}||f(t)||_\epsilon\big],
  \]
    where $C_{\rho, q, B, \epsilon}$ is a constant, depending only on $\rho, q$, $B$ and $\epsilon$.
\end{thm}
\proof
We begin the proof of  inequality \eqref{se} by establishing an estimates for $Aw(t)$, $D^{\rho}_{t}w(t)$ and then use it with equation \eqref{Cauchy3}.

To estimate $||Aw(t)||_{C(H)}$, we will use  Remark \ref{rem1} and  estimate \eqref{awt1}. Therefore we need to estimate $w(t)$ and $D_t^\rho w(t)$. Apply estimates \eqref{ES}, \eqref{RES}, Corollary \ref{cor} and \eqref{R}, \eqref{J1} to get
\[
||F(t)||_{C(H)}\leq C( ||\varphi_{1}||+||\varphi_{0}||)+C_{\rho}\bigg[\frac{\alpha^{2}|q_{k_{0}}|}{||Bq||}||\psi(t)||_{C[0,T]}
+\frac{2\alpha|q_{k_{0}}|}{||Bq||}||D^{\rho}_{t}\psi(t)||_{C[0,T]}
\]
\[
+||f{_{k_{0}}||_{C[0,T]}}  \bigg]
+C\frac{||A^{1+\epsilon}q||}{||Bq||}||\psi(t)||_{C[0,T]}+\frac{2\alpha||q||}{||Bq||}||D^{\rho}_{t}\psi(t)||_{C[0,T]}]+C\max_{0\leq t \leq T}||f(t)||_{\epsilon}\bigg].
\]
Again by estimates \eqref{RES}, \eqref{J1} we have
\[
\bigg|\bigg|\frac{1}{2}\int_{0}^{t}(t-\tau)^{\rho-1}R^{-1}\tilde{E}_{\rho,\rho}(-S^{-}(t-\tau)^{\rho})B_{1}[w(\tau)]d\tau
\]
\[
 -\frac{1}{2}\int_{0}^{t}(t-\tau)^{\rho-1}R^{-1}\tilde{E}_{\rho,\rho}(-S^{+}(t-\tau)^{\rho}B_{1}[w(\tau)]d\tau
 \]
 \[
 +\int_{0}^{t}(t-\tau)^{2\rho-1}E^{2}_{\rho,2\rho}(-\alpha(t-\tau)^{\rho})B_{1}[w_{k_0}(\tau)]d\tau v_{k_0}\bigg|\bigg|_{C(H)}
\]
\[
\leq  C_{B,q}||Aq||\int_{0}^{t}(t-\tau)^{\rho-1}\bigg(||w(\tau)||_{C(H)}+||D^{\rho}_{t}w(\tau)||_{C(H)}\bigg)d\tau
\]
\[
+C_{B,q,\rho}\bigg[|\alpha^{2}q_{k_{0}}||w_{k_{0}}(t)||_{C[0,T]}+2\alpha|q_{k_{0}}|||D^{\rho}_{t}w_{k_{0}}(t)||_{C[0,T]}\bigg].
\]
Therefore, from  \eqref{VE2} we obtain an estimate
\begin{equation}\label{esw}
||w(t)||_{C(H)}\leq C_{B,q,\epsilon}\bigg(||\varphi_{0}||+||\varphi_{1}||+||\psi(t)||_{C[0,T]}+||D^{\rho}_{t}\psi(t)||_{C[0,T]}
\end{equation}
\[
+\max_{0\leq t \leq T}||f(t)||_{\epsilon}\bigg )
+\int_{0}^{t}(t-\tau)^{\rho-1}\bigg(||w(\tau)||_{C(H)}+||D^{\rho}_{\tau}w(\tau)||_{C(H)}\bigg)d\tau.
\]
Let us now estimate $||D^{\rho}_{t}w(t)||_{C(H)}$. We will proceed as follows. Apply estimates (\ref{SES1}), \eqref{SRES}, Corollary \ref{cor} and \eqref{SR}, Lemma \ref{Int} to get
\[
||D^{\rho}_{t}w(t)||_{C(H)}\leq C_{B,q,\epsilon}\bigg(||\varphi_{0}||+||\varphi_{1}||_{\frac{1}{2}}+||\psi(t)||_{C[0,T]}+||D^{\rho}_{t}\psi(t)||_{C[0,T]}
\]
\[
+\max_{0\leq t \leq T}||f(t)||_{\epsilon}\bigg )
+\int_{0}^{t}(t-\tau)^{\rho-1}\bigg(||w(\tau)||_{C(H)}+||D^{\rho}_{\tau}w(\tau)||_{C(H)}\bigg)d\tau.
\]
Combining these estimates we have 
\[
||w(t)||_{C(H)}+||D^{\rho}_{t}w(t)||_{C(H)}\leq C_{B,q,\epsilon}\bigg(||\varphi_{0}||+||\varphi_{1}||_{\frac{1}{2}}+||\psi(t)||_{C[0,T]}+||D^{\rho}_{t}\psi(t)||_{C[0,T]}
\]
\[
+\max_{0\leq t \leq T}||f(t)||_{\epsilon}\bigg )
+\int_{0}^{t}(t-\tau)^{\rho-1}\bigg(||w(\tau)||_{C(H)}+||D^{\rho}_{\tau}w(\tau)||_{C(H)}\bigg)d\tau.
\]
for all $t\in [0, T]$. The Gronwall inequality \cite{ASH} implies
\[
||w(t)||_{C(H)}+||D^{\rho}_{t}w(t)||_{C(H)}\leq C_{B,q,\epsilon} e^{\frac{t^{\rho}}{\rho}}\bigg(||\varphi_{0}||+||\varphi_{1}||_{\frac{1}{2}}+||\psi(t)||_{C[0,T]}
\]
\[
+||D^{\rho}_{t}\psi(t)||_{C[0,T]}
+\max_{0\leq t \leq T}||f(t)||_{\epsilon}\bigg ).
\]

Finally, in order to obtain  estimate \eqref{esw}, it remains to note that   
\begin{equation*}
 (D_{t}^{\rho })^{2}w(t)=-2\alpha D_{t}^{\rho }w(t)-Aw(t)-\mu(t)Aq-2\alpha qD^{\rho}_{t}\mu(t)+f(t)
\end{equation*}
and use the estimates 
\begin{equation*}
    \lVert \mu(t)\rVert_{C[0,T]}\le C_{B,q}\left(||\psi||_{C[0,T]}+||w(t)||_{C(H)}\right),
\end{equation*}
\begin{equation*}
    \lVert D^{\rho}_{t}\mu(t)\rVert_{C[0,T]}\le C_{B,q}\left(||D^{\rho}_{t}\psi||_{C[0,T]}+||D^{\rho}_{t}w(t)||_{C(H)}\right).
\end{equation*}
Combining all the above estimates we have the assertion of the Theorem \ref{estimate3}.
\endproof
\proof of Theorem \ref{estimate}

Apply \eqref{p} to get
\[
||p(t)||_{C[0, T]}\leq C_{q, B}\big[ ||(D^{\rho}_{t})^{2} w(t)||_{C(H)}+||(D^{\rho}_{t})^{2} \psi||_{C[0, T]}\big].
\]

Equations \eqref{wq} and \eqref{mu} imply

\[
(D^{\rho}_{t})^{2}u(t)=(D^{\rho}_{t})^{2} w(t)+p(t) q.
\]
Hence, from estimates of $(D^{\rho}_{t})^{2} w(t)$ and $p(t)$, we obtain an estimate for $(D^{\rho}_t)^{2} u(t)$. On the other hand, by virtue of equation \eqref{prob1}, we will have
\[
Au(t)=-(D_{t}^{\rho})^{2} u(t)-2\alpha D_{t}^{\rho} u(t) + p(t) q +f(t).
\]

Now, to establish estimate \eqref{st}, it suffices to use the statement of Theorem \ref{estimate3}.

\endproof
\section{Conclusion} 
The fractional telegraph equation is the subject of research in a number of works by specialists. An analysis of published works shows that the inverse problem of determining the right-hand side of the equation - the source function has remained outside the attention of researchers, despite the importance of such a problem for applications. The present paper is devoted to the study of precisely this inverse problem. We also note that the additional condition for such an inverse problem is taken in the most general form, which includes almost all additional conditions involved in a similar inverse problem for subdiffusion equations. In addition, the choice of specific self-adjoint operators for $A$ in problem (\ref{prob1}) allows us to consider various mathematical models. As an example, we can take the equations of mathematical physics considered in Chap. 6 of the paper by Ruzhansky et al. \cite{Ruzhansky}, including the classical Sturm–Liouville problem, differential models with involution, fractional Sturm–Liouville operators, harmonic oscillators, Landau Hamiltonians, and fractional Laplacians. 

In almost all works where the fractional telegraph equation is studied, the highest derivative with respect to time is taken in the form $D_t^{2\rho}$. We, following the fundamental work \cite{RC}, took the highest derivative in the form of the square of the Caputo fractional derivative: $(D_t^{\rho})^2$. In this regard, we note that the cited work \cite{RC} provides a convincing motivation for such a choice. Nevertheless, the subject of our next study will be the inverse problem of determining the right-hand side of the telegraph equation, in which the derivative $D_t^{2\rho}$ participates.

\section{Acknowledgement}
The authors are grateful to A. O. Ashyralyev for posing the
problem and they convey thanks to Sh. A. Alimov for discussions of
these results.

\end{document}